\newcount\ite\ite=1\def\0{\global\ite=1\1}
\def\1{\item{\rm(\romannumeral\the\ite)}\advance\ite1\quad}
\def\phi{\varphi}

\documentclass[12pt,twoside]{article}
\usepackage{amssymb}

\font\teneufm=eufm10 scaled \magstep1
\font\seveneufm=eufm7 scaled \magstep1
\font\fiveeufm=eufm5  scaled \magstep1
\newfam\eufmfam

\textfont\eufmfam=\teneufm
\scriptfont\eufmfam=\seveneufm
\scriptscriptfont\eufmfam=\fiveeufm
\def\frak#1{{\fam\eufmfam\relax#1}}

\newfam\msbfam
\font\tenmsb=msbm10 scaled \magstep1  \textfont\msbfam=\tenmsb
\font\sevenmsb=msbm7 scaled \magstep1 \scriptfont\msbfam=\sevenmsb
\font\fivemsb=msbm5 scaled \magstep1  \scriptscriptfont\msbfam=\fivemsb

\def\dd#1{\raise1.5pt\hbox{$\,\partial\!$}/\raise-2.5pt\hbox{$\!\partial#1\,$}}

\def\tilde{\widetilde}
\def\hat{\widehat}
\def\5#1{{\mathcal #1}}

\def\AA{{\mathbb A}}
\def\RR{{\mathbb R}}
\def\CC{{\mathbb C}}

\def\TT{{\mathbb T}}
\def\PP{{\mathbb P}}

\def\FF{{\mathbb F}}
\def\GG{{\mathbb G}}

\def\ra{\rightarrow}

\def\GL{\mathop{\rm GL}\nolimits}

\def\Aut{\mathop{\rm Aut}\nolimits}

\def\Ann{\mathop{\rm Ann}\nolimits}

\def\Aff{\mathop{\rm Aff}\nolimits}
\def\L{\mathop{\rm L}\nolimits}

 \def\HollowBoxx #1#2#3{{\dimen0=#1 \advance\dimen0 by -#2
       \dimen1=#1 \advance\dimen1 by #3
        \vrule height 0pt depth #3 width #2
       \hskip -#3
       \vrule height #1 depth #3 width #3}}
 \def\LeftContraction{\mathord{\kern1.45pt \HollowBoxx{6pt}{3.5pt}{.4pt}}\,}

 \def\HollowBox #1#2#3{{\dimen0=#1 \advance\dimen0 by -#3
       \dimen1=#1 \advance\dimen1 by #3
        \vrule height #1 depth #3 width #3
        \vrule height 0pt depth #3 width #2
        \hskip -#3}}
 \def\RightContraction{\mathord{\, \HollowBox{6pt}{3.1pt}{.4pt}} \kern1.6pt}

\def\qed{{\hfill $\Box$}}
\newtheorem{theorem}{THEOREM}[section]
\newtheorem{corollary}[theorem]{Corollary}
\newtheorem{lemma}[theorem]{Lemma}

\newtheorem{remark}[theorem]{Remark}
\newtheorem{proposition}[theorem]{Proposition}

\begin{document}

\begin{center}
{\Large \bf On the Affine Homogeneity of\\
\vspace{0.3cm}
Algebraic Hypersurfaces\\
\vspace{0.4cm}
Arising from Gorenstein Algebras}
\medskip\\
\normalsize A. V. Isaev
\end{center}

\begin{quotation} 
{\small \sl \noindent 
To every Gorenstein algebra $A$ of finite vector space dimension greater than 1 over a field 
$\FF$ of characteristic zero, and a linear projection $\pi$ on its maximal ideal 
${\mathfrak m}$ with range equal to the annihilator $\Ann({\mathfrak m})$ 
of ${\mathfrak m}$, one can associate a certain algebraic hypersurface $S_{\pi}\subset{\mathfrak m}$. Such hypersurfaces possess remarkable properties. They can be used, for instance, to help decide whether two given Gorenstein algebras are isomorphic, which for $\FF=\CC$ leads to interesting consequences in singularity theory. Also, for $\FF=\RR$ such hypersurfaces naturally arise in CR-geometry. Applications of these hypersurfaces to problems in algebra and geometry are particularly striking when the hypersurfaces are affine homogeneous. In the present paper we establish a criterion for the affine homogeneity of $S_{\pi}$. This criterion requires the automorphism group $\Aut({\mathfrak m})$ of ${\mathfrak m}$ to act transitively on the set of hyperplanes in ${\mathfrak m}$ complementary to $\Ann({\mathfrak m})$. As a consequence of this result we obtain the affine homogeneity of $S_{\pi}$ under the assumption that the algebra $A$ is graded.}
\end{quotation}

\thispagestyle{empty}

\pagestyle{myheadings}
\markboth{A. V. Isaev}{Algebraic Hypersurfaces Arising from Gorenstein Algebras}

\setcounter{section}{0}

\section{Introduction}\label{intro}
\setcounter{equation}{0}

In this paper we consider finite-dimensional Gorenstein algebras over a field $\FF$ of characteristic zero. Recall that a local commutative associative algebra $A$ of finite vector space dimension greater than 1 is Gorenstein if and only if the annihilator $\Ann({\mathfrak m})$ of its maximal ideal ${\mathfrak m}$ is 1-dimensional (see, e.g. \cite{Hu}). In the case $\FF=\CC$, in our earlier paper \cite{FIKK} we found a criterion for two such algebras to be isomorphic. The criterion was given in terms of a certain algebraic hypersurface $S_{\pi}\subset{\mathfrak m}$ associated to a linear projection $\pi$ on ${\mathfrak m}$ with range $\Ann({\mathfrak m})$. The hypersurface $S_{\pi}$ passes through the origin and is the graph of a polynomial map $P_{\pi}:\ker\pi\ra \Ann({\mathfrak m})\simeq\CC$ that has no linear term and is completely determined by its second- and third-order terms. In Proposition 2.3 of \cite{FIKK} we showed that for $\FF=\CC$ two Gorenstein algebras $A$, $\tilde A$ are isomorphic if and only if any two hypersurfaces $S_{\pi}$ and $S_{\tilde\pi}$ arising from $A$ and $\tilde A$, respectively, are affinely equivalent, that is, if there exists a bijective affine map ${\mathcal A}: {\mathfrak m}\ra\tilde {\mathfrak m}$ such that ${\mathcal A}(S_{\pi})=S_{\tilde\pi}$. The proof is based on a CR-geometric argument and also works for $\FF=\RR$. 

As we will see in Section \ref{prelim}, the hypersurface $S_{\pi}$ can be introduced for Gorenstein algebras over any field $\FF$ of characteristic zero. The first result of this paper is an extension of Proposition 2.3 of \cite{FIKK} to the case of an arbitrary field of zero characteristic (see Proposition \ref{equivalence} in Section \ref{prelim}). Our proof is completely algebraic and requires no CR-geometry even in the cases $\FF=\RR$, $\FF=\CC$. The criterion for $A$ and $\tilde A$ to be isomorphic given in Proposition \ref{equivalence} takes a particularly nice form if at least one of $S_{\pi}$, $S_{\tilde\pi}$ is {\it affine homogeneous} (recall that a subset ${\mathcal S}$ of an affine space $V$ is called affine homogeneous if for every pair of points $p,q\in{\mathcal S}$ there exists a bijective affine map ${\mathcal A}$ of $V$ such that ${\mathcal A}({\mathcal S})={\mathcal S}$ and ${\mathcal A}(p)=q$). In this case the hypersurfaces $S_{\pi}$, $S_{\tilde\pi}$ are affinely equivalent if and only if they are linearly equivalent, that is, equivalent by means of a bijective linear transformation. As explained in \cite{FIKK}, for any field $\FF$ the linear equivalence of $S_{\pi}$, $S_{\tilde\pi}$ takes place if and only if the corresponding polynomials $P_{\pi}$, $P_{\tilde\pi}$ are linearly equivalent up to scale. To establish the equivalence or non-equivalence of $P_{\pi}$, $P_{\tilde\pi}$ it is sufficient to consider their quadratic and cubic terms alone. Thus, if at least one of $S_{\pi}$, $S_{\tilde\pi}$ is affine homogeneous, verifying that $A$ and $\tilde A$ are isomorphic reduces to verifying that certain  pairs of quadratic and cubic forms are linearly equivalent up to scale. For the case $\FF=\CC$ these results were summarized in Theorem 2.11 of \cite{FIKK}, which we subsequently applied to the moduli algebras of complex quasi-homogeneous isolated hypersurface singularities. We then used the Mather-Yau theorem (see \cite{MY}) to deduce equivalence results for such singularities. For many examples solving the equivalence problem for the polynomials $P_{\pi}$, $P_{\tilde\pi}$ has turned out to be much easier than solving the equivalence problem for the algebras $A$, $\tilde A$ directly (see, e.g. \cite{I2} for the case of simple elliptic singularities).  

As stated above, the affine homogeneity property of $S_{\pi}$ played an important role in our arguments in \cite{FIKK}. In Proposition 2.4 of \cite{FIKK} this property was established in the case of graded Gorenstein algebras over $\CC$, which was sufficient for our application to isolated hypersurface singularities. In this paper, we investigate the affine homogeneity of $S_{\pi}$ for an arbitrary field $\FF$ of characteristic zero. First, we establish a criterion for the affine homogeneity of $S_{\pi}$ in Theorem \ref{main1}. Interestingly, it turns out that $S_{\pi}$ is affine homogeneous if and only if the automorphism group $\Aut({\mathfrak m})$ of the algebra ${\mathfrak m}$ acts transitively on the set $\TT$ of all hyperplanes in ${\mathfrak m}$ complementary to $\Ann({\mathfrak m})$, in which case we say that the algebra $A$ has {\it Property} (P). 

The above result indicates the importance of studying the group $\Aut({\mathfrak m})$, which is an affine algebraic group over $\FF$. Observe that, as an algebraic variety, $\TT$ is isomorphic to the affine space $\AA^n$ over $\FF$, with $n:=\dim {\mathfrak m}-1$, and therefore Property (P) implies the lower dimension bound $\dim\Aut({\mathfrak m})\ge n$. This estimate is known to hold for any graded Gorenstein algebra (see \cite{XY}). The only lower bound known for general Gorenstein algebras is $\dim\Aut({\mathfrak m})\ge\dim({\mathfrak m}/{\mathfrak m}^2)$ (see, e.g. \cite{Pe}).

Next, we obtain Property (P) for any graded Gorenstein algebra over a field of zero characteristic (see Theorem \ref{main2}). It then follows that the affine homogeneity of $S_{\pi}$ takes place for any such algebra (see Corollary \ref{cor2}). We stress that the proof of the affine homogeneity of $S_{\pi}$ given in this paper is not based on the proof presented in \cite{FIKK} for $\FF=\CC$.

The hypersurface $S_{\pi}$ is a rather fascinating object from both the algebraic and geometric viewpoints. Indeed, as we explained above, such hypersurfaces can be used to establish the equivalence or non-equivalence of Gorenstein algebras. Next, they are of interest to affine geometry since, at least in the graded case, they represent a way to associate an affine homogeneous affine algebraic variety to any Gorenstein algebra. Furthermore, for $\FF=\RR$ these hypersurfaces arise in CR-geometry as the bases of spherical tube hypersurfaces. In fact, $S_{\pi}$ was first introduced in \cite{FK1} for the purpose of classifying spherical tube hypersurfaces up to affine equivalence. The literature on spherical tube hypersurfaces is quite extensive, and more information on this subject, including numerous examples where $S_{\pi}$ was explicitly computed, can be found, for instance, in our articles \cite{I1}, \cite{I3} and forthcoming monograph \cite{I4}. We note that in \cite{I3} we obtained the affine homogeneity of $S_{\pi}$ for some special graded Gorenstein algebras over $\RR$, but the proof in fact works for an arbitrary field.

In the next section we give all necessary definitions and state our main results. Proofs are presented in Sections \ref{proofs1} and \ref{proofs2}. 

{\bf Acknowledgements.} This work is supported by the Australian Research Council. Similar results were independently obtained in recent preprint \cite{FK2} by different methods. We would like thank Nikolay Kruzhilin for useful discussions.

\section{Statement of Results}\label{prelim}
\setcounter{equation}{0}
    
Let $A$ be a Gorenstein algebra of finite vector space dimension greater than 1 over a field $\FF$ of characteristic zero, ${\mathfrak m}$ the maximal ideal of $A$, $n:=\dim{\mathfrak m}-1$ (observe that $n\ge 0$), and $\pi$ a linear projection on ${\mathfrak m}$ with range $\Ann({\mathfrak m}):=\{u\in{\mathfrak m}: u\cdot{\mathfrak m}=0\}$ (we call such projections {\it admissible}). Consider the following $\Ann({\mathfrak m})$-valued bilinear form on $A$:
\begin{equation}
b_{\pi}(a,c):=\pi(ac), \quad a,c\in A,\label{b}
\end{equation}
where $\pi$ is extended to all of $A$ by the condition $\pi({\bf 1})=0$, with ${\bf 1}$ being the unit of $A$. It is well-known that the form $b_{\pi}$ is non-degenerate (see, e.g. \cite{He}, p. 11). Further, let $\PP(A)$ be the projectivization of $A$, and consider the following projective quadric:
$$
Q_{\pi}:=\Bigl\{[a]\in\PP(A):b_{\pi}(a,a)=0\Bigr\},
$$
where $[a]$ denotes the point of $\PP(A)$ represented by $a\in A$. The inclusion ${\mathfrak m}\subset A$ induces the inclusion $\PP({\mathfrak m})\subset \PP(A)$, and we think of $\PP({\mathfrak m})$ as the hyperplane at infinity in $\PP(A)$. We then identify ${\bf 1}+{\mathfrak m}\subset A$ with the affine part of $\PP(A)$ and introduce the corresponding affine quadric $Q'_{\pi}:=Q_{\pi}\cap ({\bf 1}+{\mathfrak m})$. Below we often, without stating it explicitly, identify ${\bf 1}+{\mathfrak m}$ with ${\mathfrak m}$ in the obvious way, and upon this identification one has
\begin{equation}
Q_{\pi}'=\{u\in{\mathfrak m}: \pi(2u+u^2)=0\}=\{u\in{\mathfrak m}: 2u+u^2\in{\mathcal K}\},\label{equation}
\end{equation}
where ${\mathcal K}:=\ker\pi$.

Let $\exp: {\mathfrak m}\ra {\bf 1}+{\mathfrak m}$ be the exponential map
$$
\displaystyle\exp(u):={\bf 1}+\sum_{m=1}^{\infty}\frac{1}{m!}u^m.
$$
By Nakayama's lemma the maximal ideal ${\mathfrak m}$ of $A$ is a nilpotent algebra, and therefore the above sum is in fact finite, with the highest-order term corresponding to $m=\nu$, where $\nu\ge 1$ is the nil-index of ${\mathfrak m}$ (i.e. the largest integer $\ell$ for which ${\frak m}^{\ell}\ne 0$). Thus, the exponential map is a polynomial transformation. It is bijective with the inverse given by the polynomial transformation
$$
\log({\bf 1}+u):=\sum_{m=1}^{\infty}\frac{(-1)^{m+1}}{m}u^m,\quad u\in{\mathfrak m}.
$$

We are now ready to introduce the main object of our study. Define
$$
S_{\pi}:=\log(Q_{\pi}')\subset{\mathfrak m}. 
$$
Observe that $S_{\pi}$ is an algebraic hypersurface in ${\mathfrak m}$ passing through the origin. It is easily seen that $S_{\pi}$ is the graph of a polynomial map $P_{\pi}:{\mathcal K}\ra\Ann({\mathfrak m})$. 

Clearly, the hypersurface
$S_{\pi}$ depends on the choice of the admissible projection $\pi$. It is easy to observe (using Lemma \ref{twopi} below) that for two admissible projections $\pi$ and $\pi'$ on ${\mathfrak m}$ one has
\begin{equation}
S_{\tilde\pi}=S_{\pi'}+x\label{translation}
\end{equation}
for some $x\in{\mathfrak m}$. 

As stated in the following proposition, the isomorphism class of the algebra $A$ is completely determined by the affine equivalence class of $S_{\pi}$.  

\begin{proposition}\label{equivalence}\sl Let $A$, $\tilde A$ be Gorenstein algebras of finite vector space dimension greater than 1, and $\pi$, $\tilde\pi$ admissible projections on $A$, $\tilde A$, respectively. Then $A$, $\tilde A$ are isomorphic if and only if the hypersurfaces $S_{\pi}$, $S_{\tilde\pi}$ are affinely equivalent. Moreover, if ${\mathcal A}:{\mathfrak m}\ra\tilde{\mathfrak m}$ is an affine map such that ${\mathcal A}(S_{\pi})=S_{\tilde\pi}$, then the linear part of ${\mathcal A}$ is an algebra isomorphism between ${\mathfrak m}$ and $\tilde{\mathfrak m}$.
\end{proposition}

\noindent {\bf Proof:} If $A$ and $\tilde A$ are isomorphic, then $S_{\pi}$ and $S_{\tilde\pi}$ are affinely equivalent by (\ref{translation}). Conversely, let ${\mathcal A}:{\mathfrak m}\ra\tilde{\mathfrak m}$ be an affine equivalence with ${\mathcal A}(S_{\pi})=S_{\tilde\pi}$, and $y:={\mathcal A}(0)$. Then the linear map ${\mathcal L}(u):={\mathcal A}(u)-y$, with $u\in{\mathfrak m}$, maps $S_{\pi}$ onto $S_{\tilde\pi}-y$. Consider the admissible projection on $\tilde{\mathfrak m}$ given by $\tilde\pi'(v):=\tilde\pi(\widetilde{\exp}(2y)v)$, where $\widetilde{\exp}$ is the exponential map associated to $\tilde A$. We then have $
S_{\tilde\pi}-y=S_{\tilde\pi'}$, hence ${\mathcal L}$ maps $S_{\pi}$ onto $S_{\tilde\pi'}$. 

Recall that $S_{\pi}$ is the graph of a polynomial map $P_{\pi}:{\mathcal K}\ra\Ann({\mathfrak m})$. Choose coordinates $\alpha=(\alpha_1,\dots,\alpha_n)$ in ${\mathcal K}$ and a coordinate $\alpha_0$ in $\Ann({\mathfrak m})$. In these coordinates the hypersurface $S_{\pi}$ is written as
$$
\alpha_0=P_{\pi}(\alpha_1,\dots,\alpha_n)=\sum_{i,j=1}^ng_{ij}\alpha_i\alpha_j+\sum_{i,j,k=1}^nh_{ijk}\alpha_i\alpha_j\alpha_k+\dots,
$$ 
where $g_{ij}$ and $h_{ijk}$ are symmetric in all indices, $(g_{ij})$ is non-degenerate, and the dots denote higher-order terms. In Proposition 2.10 in \cite{FIKK} we showed that the above equation of $S_{\pi}$ is in Blaschke normal form, that is, one has $\sum_{ij}g^{ij}h_{ijk}=0$ for all $k$, where $(g^{ij}):=(g_{ij})^{-1}$ (observe that the proof of this statement in \cite{FIKK} works for any field of characteristic zero). As shown in Proposition 1 of \cite{EE} (which is valid over any field of characteristic zero), every linear isomorphism between hypersurfaces in Blaschke normal form is written as
$$
\alpha_0\mapsto\hbox{const}\,\alpha_0,\quad \alpha\mapsto C\alpha,
$$
where $C\in\GL(n,\FF)$. Hence, representing $S_{\tilde\pi}$ in coordinates $\tilde\alpha_0$,\linebreak $\tilde\alpha=(\tilde\alpha_1,\dots,\tilde\alpha_n)$ in $\tilde{\mathfrak m}$ chosen as above, we see that the map  ${\mathcal L}$ has the form
$$
\tilde\alpha_0=\hbox{const}\,\alpha_0,\quad \tilde\alpha= C\alpha
$$
for some $C\in\GL(n,\FF)$. Therefore, the polynomials $P_{\pi}$ and $P_{\tilde\pi}$ are linearly equivalent up to scale. 

As explained in Lemma 9.30 in \cite{FK1}, the algebras ${\mathfrak m}$ and $\tilde{\mathfrak m}$ can be completely reconstructed from the second- and third-order terms of the polynomials $P_{\pi}$ and $P_{\tilde\pi}$, respectively. Since these polynomials are linearly equivalent up to scale, the algebras ${\mathfrak m}$ and $\tilde{\mathfrak m}$ are isomorphic, and, moreover, the linear map ${\mathcal L}$ is an algebra isomorphism between them (see Proposition 9.33 in \cite{FK1}). The proof is complete.\qed    
\vspace{0.1cm}\\

We will now concentrate on the affine homogeneity property of $S_{\pi}$. Since (\ref{translation}) implies that $S_{\pi'}$ is affine homogeneous if and only if $S_{\pi}$ is affine homogeneous, the projection $\pi$ will be fixed from now on.

Let $\Aff({\mathfrak m})$ be the group of all bijective affine transformations of ${\mathfrak m}$. For any subset ${\mathcal S}\subset{\mathfrak m}$ define
$$
\Aff({\mathcal S}):=\left\{f\in\Aff({\mathfrak m}): f({\mathcal S})={\mathcal S}\right\}.
$$
The hypersurface $S_{\pi}$ is affine homogeneous if and only if the group $\Aff(S_{\pi})$ acts on $S_{\pi}$ transitively. Further, let $\GG$ be the Grassmannian of all hyperplanes in ${\mathfrak m}$, and $\TT$ the subset of $\GG$ that consists of all hyperplanes in ${\mathfrak m}$ complementary to $\Ann({\mathfrak m})$. As an algebraic variety, $\TT$ is isomorphic to the affine space $\AA^n$ over $\FF$. Clearly, $\Aut({\mathfrak m})$ acts on $\TT$, and the action is algebraic. We say that the algebra $A$ has Property (P) if the action of $\Aut({\mathfrak m})$ on $\TT$ is transitive.

We are now ready to state our first main theorem.   

\begin{theorem}\label{main1}\sl The group $\Aff(S_{\pi})$ acts transitively on $S_{\pi}$ if and only if the algebra $A$ has Property {\rm (P)}. 
\end{theorem}

\noindent Theorem \ref{main1} will be proved in two steps. At the first step we introduce the subgroup
\begin{equation}
G_{\pi}:=\{f\in\Aff(S_{\pi}):f=\log\circ\, g\,\circ\exp\,\hbox{for some}\, g\in\Aff(Q_{\pi}')\},\label{specgroup}
\end{equation}
and prove Theorem \ref{main1} for $G_{\pi}$ in place of  $\Aff(S_{\pi})$ (see Proposition \ref{prop1}). At the second step we show that $G_{\pi}=\Aff(S_{\pi})$ (see Proposition \ref{prop2}).

\begin{remark}\label{rem} \rm Observe that the affine quadric $Q_{\pi}'$ is affine homogeneous. Indeed, let $\alpha=(\alpha_1,\dots,\alpha_n)$ be coordinates in ${\mathcal K}$ and $\alpha_0$ a coordinate in $\Ann({\mathfrak m})$. Then $Q_{\pi}'$ is given by the equation
$$
\alpha_0=\langle \alpha,\alpha\rangle,
$$
where $\langle \cdot, \cdot\rangle$ is a non-degenerate quadratic form (see (\ref{equation})), and the group of affine transformations of ${\mathfrak m}$
$$
\begin{array}{l}
\alpha_0^*=\alpha_0+2\langle \alpha,\beta\rangle + \langle \beta,\beta\rangle,\\
\vspace{-0.1cm}\\
\alpha^*=\alpha+\beta,\quad \beta\in\FF^n,
\end{array}
$$
acts transitively on $Q_{\pi}'$. For every $g\in\Aff(Q_{\pi}')$ the composition $f:=\log\circ\, g\circ\exp$ is a polynomial map of ${\mathfrak m}$ having a polynomial inverse, and thus $S_{\pi}$ is homogeneous under such maps. However, it is easy to find an example of $g\in\Aff(Q_{\pi}')$ for which the corresponding map $f$ is not affine.
\end{remark} 

Theorem \ref{main1} is proved in Section \ref{proofs1}. The proof also shows that the groups $\Aff(S_{\pi})$ and $\Aut({\mathfrak m})$ are canonically isomorphic. Furthermore, the isomorphism between $\Aff(S_{\pi})$ and $\Aut({\mathfrak m})$ identifies the action of $\Aff(S_{\pi})$ on $S_{\pi}$ with the action of $\Aut({\mathfrak m})$ on $\TT$ (see Remark \ref{isomorphism}).

Next, we say that a Gorenstein algebra $A$ of finite vector space dimension greater than 1 is graded if
\begin{equation}
A=\bigoplus_{j\ge0}A_{j},\quad A_{j}A_{k}\subset A_{j+k},\label{grading}
\end{equation}
where $A_{j}$ are linear subspaces of $A$, with $A_0\simeq\FF$. Then ${\mathfrak m}=\oplus_{j>0}A_j$ and $\Ann({\mathfrak m})=A_d$ for $d:=\max\{j:A_{j}\ne 0\}$. 

We will now state our second main result.

\begin{theorem}\label{main2}\sl Every graded Gorenstein algebra has Property {\rm (P)}.
\end{theorem}

\noindent We prove Theorem \ref{main2} in Section \ref{proofs2}. Clearly, together with Theorem \ref{main1} it yields the following corollary.

\begin{corollary}\label{cor2}\sl If $A$ is graded, then $S_{\pi}$ is affine homogeneous.
\end{corollary}

\noindent We will now give proofs of the two theorems stated in this section.  
   
\section{Proof of Theorem \ref{main1}}\label{proofs1}
\setcounter{equation}{0}

We will first prove the following proposition.

\begin{proposition}\label{prop1}\sl Let $G_{\pi}\subset\Aff(S_{\pi})$ be the subgroup defined in {\rm (\ref{specgroup})}. The subgroup $G_{\pi}$ acts transitively on $S_{\pi}$ if and only if the algebra $A$ has\linebreak Property {\rm (P)}.
\end{proposition}

We start with two lemmas. Let $g\in\Aff({\mathfrak m})$ and $f:=\log\circ\, g\circ\exp$. Clearly, $f$ is a polynomial transformation of ${\mathfrak m}$ having a polynomial inverse. In the first lemma we find a necessary and sufficient condition for $f$ to be an affine map. Write $g=L_g+u_g$, where $L_g$ is the linear part of $g$ and $u_g:=g(0)$.

\begin{lemma}\label{baffine}\sl The map $f$ is affine if and only if the linear map of ${\mathfrak m}$
\begin{equation}
L: u\mapsto ({\bf 1}+u_g)^{-1}\cdot L_g(u)\label{hom}
\end{equation}
is an automorphism of the algebra ${\mathfrak m}$. In this case we have $f=L+\log({\bf 1}+u_g)$.
\end{lemma}

\noindent{\bf Proof:} Suppose first that $L\in\Aut({\mathfrak m})$. To show that $f=L+\log({\bf 1}+u_g)$, we need to check that $\exp\circ \bigl(L+\log({\bf 1}+u_g)\bigr)=g\circ \exp$, where, as usual, ${\bf 1}+{\mathfrak m}$ is identified with ${\mathfrak m}$. We have
$$
\begin{array}{l}
\displaystyle\exp\circ \Bigl(L(u)+\log({\bf 1}+u_g)\Bigr)=({\bf 1}+u_g)\exp(L(u))=\\
\vspace{-0.1cm}\\
\displaystyle({\bf 1}+u_g)\left({\bf 1}+\sum_{m=1}^{\infty}\frac{(L(u))^m}{m!}\right)=({\bf 1}+u_g)\left({\bf 1}+\sum_{m=1}^{\infty}\frac{L(u^m)}{m!}\right)=\\
\vspace{-0.1cm}\\
\displaystyle {\bf 1}+u_g+\sum_{m=1}^{\infty}\frac{L_g(u^m)}{m!}={\bf 1}+L_g\left(\sum_{m=1}^{\infty}\frac{u^m}{m!}\right)+u_g=g\circ \exp,
\end{array}
$$
as required.

Suppose now that $f$ is affine and write it as $L_f+u_f$, where $L_f$ is the linear part of $f$ and $u_f:=f(0)$. We have
\begin{equation}
\exp\circ f=g\circ \exp.\label{imp}
\end{equation}
Comparing the constant terms in identity (\ref{imp}), we see $u_f=\log({\bf 1}+u_g)$. Next, by comparing the first-order terms we obtain $L_f=L$. Further, if $n\ge 1$ then $\nu\ge 2$ (recall that $\nu$ is the nil-index of ${\mathfrak m}$), and comparison of the second-order terms yields
\begin{equation}
L(u^2)\equiv (L(u))^2.\label{square}
\end{equation}
If $n=0$ then ${\mathfrak m}=\Ann({\mathfrak m})$, and (\ref{square}) is obvious. Clearly, identity (\ref{square}) implies $L\in\Aut({\mathfrak m})$. Observe also that for $\nu\ge 3$ comparison of higher-order terms in (\ref{imp}) provides no additional information about $f$. The proof is complete.\qed
\vspace{0.1cm}\\

In the second lemma we find a necessary and sufficient condition for $g$ to be an element of $\Aff(Q_{\pi}')$ under the assumptions that $f$ is affine and $u_g\in Q_{\pi}'$.
  
\begin{lemma}\label{selfmap} \sl Let $f$ be affine and $u_g\in Q_{\pi}'$. Then $g\in\Aff(Q_{\pi}')$ if and only if
\begin{equation}
({\bf 1}+u_g)^2L({\mathcal K})={\mathcal K},\label{idensss}
\end{equation}
where $L$ is defined in {\rm (\ref{hom})}.
\end{lemma}

\noindent {\bf Proof:} According to (\ref{equation}), we have $g(Q_{\pi}')\subset Q_{\pi}'$ if and only if 
\begin{equation}
2(L_g(u)+u_g)+(L_g(u)+u_g)^2\in{\mathcal K}\,\,\hbox{whenever $2u+u^2\in{\mathcal K}$.}\label{idens}
\end{equation}
Since $u_g\in Q_{\pi}'$, $L_g=({\bf 1}+u_g)\cdot L$, and $L\in\Aut({\mathfrak m})$, condition (\ref{idens}) is equivalent to
\begin{equation}
({\bf 1}+u_g)^2L(2u+u^2)\in{\mathcal K}\,\,\hbox{whenever $2u+u^2\in{\mathcal K}$.}\label{idens1}
\end{equation}
Observe now that every element $v\in{\mathfrak m}$ can be represented in the form $v=2u+u^2$ for some $u\in{\mathfrak m}$. Indeed, take $u:=({\bf 1}+v)^{1/2}-{\bf 1}$ where
$$
\displaystyle ({\bf 1}+v)^{1/2}:={\bf 1}+\sum_{m=1}^{\infty}\frac{(-1)^{m+1}\cdot 3 \cdot 5 \cdot \dots \cdot (2m-3)}{2^m\, m!}v^m.
$$
Therefore, condition (\ref{idens1}) is equivalent to identity (\ref{idensss}). Thus, we have shown that $g(Q_{\pi}')\subset Q_{\pi}'$ if and only if (\ref{idensss}) holds.

To complete the proof of the lemma we need to show that identity (\ref{idensss}) implies $g^{-1}(Q_{\pi}')\subset Q_{\pi}'$. By (\ref{equation}), the inclusion  $g^{-1}(Q_{\pi}')\subset Q_{\pi}'$ holds if and only if 
$$
2(L_g^{-1}(u)-L_g^{-1}(u_g))+(L_g^{-1}(u)-L_g^{-1}(u_g))^2\in{\mathcal K}\,\,\hbox{whenever $2u+u^2\in{\mathcal K}$.}
$$
This condition follows by observing that $L_g^{-1}=({\bf 1}+u_g)^{-1}L^{-1}$ and\linebreak $({\bf 1}+u_g)^{-2}L^{-1}({\mathcal K})={\mathcal K}$.\qed
\vspace{0.1cm}\\

\noindent {\bf Proof of Proposition \ref{prop1}:} The group $G_{\pi}$ acts transitively on $S_{\pi}$ if and only if the group $H_{\pi}:=\{\exp\circ f\circ\log: f\in G_{\pi}\}\subset\Aff(Q_{\pi}')$ acts transitively on $Q_{\pi}'$. Lemmas \ref{baffine}, \ref{selfmap} imply that the group $H_{\pi}$ acts transitively on $Q_{\pi}'$ if and only if for every $u\in Q_{\pi}'$ there exists $\varphi\in\Aut({\mathfrak m})$ satisfying
\begin{equation}
\varphi({\mathcal K})=({\bf 1}+(2u+u^2))^{-1}{\mathcal K}.\label{idens8}
\end{equation}
As the element $u$ varies over all points of $Q_{\pi}'$, the element $2u+u^2$ varies over all points of ${\mathcal K}$.

We now need the following lemma.

\begin{lemma}\label{twopi}\sl Let $\Pi,\tilde\Pi\in\TT$. Then there exists $x,y\in\Pi$ such that $\tilde\Pi=({\bf 1}+x)^{-1}\Pi=({\bf 1}+y)\Pi$.
\end{lemma}

\noindent {\bf Proof:} The proof is similar to that of Proposition 2.2 of \cite{FIKK}. First, observe that $\Pi=\{u+\mu(u):u\in
\tilde\Pi\}$ for some linear map $\mu:\tilde\Pi\ra\Ann({\mathfrak m})$. Clearly, the space $\L\bigl(\tilde\Pi, \Ann({\mathfrak m})\bigr)$ of all linear maps from $\tilde\Pi$ to $\Ann({\mathfrak m})$ has dimension $n=\dim\Pi=\dim\tilde\Pi$. 

Consider the linear map
$$
\Phi:\Pi\ra \L\bigl(\tilde\Pi, \Ann({\mathfrak m})\bigr),\quad\Phi(x)(u)=\pi(x u),\quad x\in\Pi,\,\,u\in\tilde\Pi,
$$
where $\pi$ is the admissible projection on ${\mathfrak m}$ with kernel $\Pi$. Since the form $b_{\pi}$ defined in (\ref{b}) is non-degenerate on $A$, the map $\Phi$ is injective and hence is an isomorphism. Therefore, $\mu(u)\equiv\pi(x u)$ for some $x\in\Pi$, and thus $\tilde\Pi=({\bf 1}+x)^{-1}\Pi$. We now write $({\bf 1}+x)^{-1}$ as ${\bf 1}+x'$ with $x':=\sum_{m=1}^{\infty}(-1)^m x^m$ and take $y$ to be the projection of $x'$ to $\Pi$ along $\Ann({\mathfrak m})$. This completes the proof of the lemma.\qed
\vspace{0.1cm}\\

Lemma \ref{twopi} shows that as the element $u$ varies over all points of $Q_{\pi}'$, the right-hand side of (\ref{idens8}) varies over all points of $\TT$. Thus, $G_{\pi}$ acts transitively on $S_{\pi}$ if and only if the algebra $A$ has Property (P). The proof of Proposition \ref{prop1} is complete. \qed
\vspace{0.1cm}\\

To finalize the proof of Theorem \ref{main1} we need the following proposition.

\begin{proposition}\label{prop2}\sl We have $G_{\pi}=\Aff(S_{\pi})$.
\end{proposition}

\noindent {\bf Proof:} We will show that $\Aff(S_{\pi})\subset G_{\pi}$. Let ${\mathcal A}$ be an element of $\Aff(S_{\pi})$ and ${\mathcal L}$ the linear part of ${\mathcal A}$. By Proposition \ref{equivalence} we have ${\mathcal L}\in\Aut({\mathfrak m})$. Let ${\mathcal K}':={\mathcal L}({\mathcal K})$. By Lemma \ref{twopi} there exists $x\in{\mathcal K}$ such that ${\mathcal K}'=({\bf 1}+x)^{-1}{\mathcal K}$. The element $x$ can be written as $x=2u+u^2$ for some $u\in{\mathcal Q}_{\pi}'$. Hence, by Lemmas \ref{baffine}, \ref{selfmap} there exists $f\in G_{\pi}$ with $L_f={\mathcal L}$. We will now show that $f={\mathcal A}$.

Suppose that $f\ne{\mathcal A}$. Then there exists a non-zero $y\in S_{\pi}$ such that the map $\hat f(u):=u+y$, with $u\in{\mathfrak m}$, lies in $\Aff(S_{\pi})$. Clearly, $\hat f=\log\circ\, \hat g\circ\exp$, where $\hat g(u):=\exp(y)u+(\exp(y)-{\bf 1})$, with $u_{\hat g}=\exp(y)-{\bf 1}\in Q_{\pi}'$. Since $\hat f\in\Aff(S_{\pi})$, we have $\hat g\in\Aff(Q_{\pi}')$. By Lemma \ref{selfmap} we then obtain $({\bf 1}+z){\mathcal K}={\mathcal K}$, where
$z:=2u_{\hat g}+u_{\hat g}^2\in{\mathcal K}$. Thus, $zu\in{\mathcal K}$ for all $u\in{\mathcal K}$, which implies $b_{\pi}(z,A)=0$. The non-degeneracy of $b_{\pi}$ now yields $z=0$ and therefore $\exp(2y)={\bf 1}$, which implies $y=0$. This contradiction concludes the proof.\qed
\vspace{0.1cm}\\

Theorem \ref{main1} now follows from Propositions \ref{prop1} and \ref{prop2}.\qed
\vspace{0.1cm}\\

\begin{remark}\label{isomorphism}\rm Our proof of Theorem \ref{main1} yields an alternative proof of Proposition 4.9 in \cite{FK2}. Indeed, let $\Psi:\Aff(S_{\pi})\ra\Aut({\mathfrak m})$ be the map that assigns every element of $\Aff(S_{\pi})$ its linear part. The above argument shows that $\Psi$ is in fact an isomorphism and the diagram
$$
\begin{array}{ccccc}
\Aut({\mathfrak m})&\hspace{-0.2cm}\times&\hspace{-0.4cm}\TT&\hspace{-0.4cm}\longrightarrow&\hspace{-0.4cm}\TT\\
\vspace{-0.1cm}\\
\hspace{0.2cm}\Big\uparrow\hbox{\small$\Psi$}&&\hspace{-0.2cm}\Big\uparrow\hbox{\small$\psi$}&&\hspace{-0.2cm}\Big\uparrow\hbox{\small$\psi$}\\
\vspace{-0.1cm}\\
\Aff(S_{\pi})&\hspace{-0.2cm}\times&\hspace{-0.4cm}S_{\pi}&\hspace{-0.4cm}\longrightarrow&\hspace{-0.4cm}S_{\pi}\\
\end{array}
$$
is commutative. Here $\psi$ is the bijective map defined by $\psi(u):=\exp(-2u){\mathcal K}$ and the horizontal arrows refer to the actions of the groups $\Aut({\mathfrak m})$ and $\Aff(S_{\pi})$ on $\TT$ and $S_{\pi}$, respectively. 
\end{remark}

\section{Proof of Theorem \ref{main2}}\label{proofs2}
\setcounter{equation}{0}

Let $A$ be a Gorenstein algebra of finite vector space dimension greater than 1, and assume in addition that $A$ is graded (see (\ref{grading})). In this case every element $u\in{\mathfrak m}$ can be uniquely written as $u=u_1+\dots+u_d$, where $u_j\in A_j$. Let $\Pi:=\oplus_{j=1}^{d-1}A_j$. Clearly, $\Pi$ lies in $\TT$. Fix another element $\tilde\Pi\in\TT$. To prove the theorem, we need to show that there exists $\varphi\in\Aut({\mathfrak m})$ such that $\varphi(\Pi)=\tilde\Pi$. By Lemma \ref{twopi} one can find $y\in\Pi$ with $\tilde\Pi=({\bf 1}+y)\Pi$. Therefore, we need to show that for any $y\in\Pi$ there exists $\varphi\in\Aut({\mathfrak m})$ for which $\varphi(\Pi)=({\bf 1}+y)\Pi$. Equivalently, $\varphi$ must be chosen so that for every $u\in\Pi$ one can find $v\in\Pi$ with $\varphi(u)=({\bf 1}+y)v$. 

We construct $\varphi$ using derivations of ${\mathfrak m}$ (cf. the proof of Proposition 2.3 of \cite{XY}). For every $x\in\Pi$ set
$$
D_{x}(u):=x(u_1+2u_2+\dots+(d-1) u_{d-1}), \quad u\in{\mathfrak m}.
$$
It is straightforward to check that $D_{x}$ is a derivation of ${\mathfrak m}$. Let $\varphi_{x}$ be the automorphism of ${\mathfrak m}$ obtained by exponentiating $D_{x}$. If $x=x_i\in A_i$ for $1\le i\le d-1$, the automorphism $\varphi_{x}$ can be easily computed:
$$
\begin{array}{l}
\displaystyle\varphi_{x_i}(u)=\sum_{j=1}^{d-1}\left({\bf 1}+jx_i+\frac{j(j+i)}{2!}x_i^2+\frac{j(j+i)(j+2i)}{3!}x_i^3+\dots\right)u_j+u_d=\\
\vspace{-0.1cm}\\
\hspace{1.68cm}\displaystyle\sum_{j=1}^{d-1}({\bf 1}-ix_i)^{-j/i}u_j+u_d, \quad u\in{\mathfrak m}.
\end{array}
$$
We look for $\varphi$ in the form $\varphi_{x_{d-1}}\circ\dots\circ\varphi_{x_1}$ and show that one can choose suitable $x_i$, $i=1,\dots,d-1$.

Fix $u\in\Pi$ and let $w:=\varphi_{x_{d-1}}\circ\dots\circ\varphi_{x_1}(u)$. It is then straightforward to see that
\begin{equation}
\hspace{-0.2cm}\begin{array}{l}
w_j=\mu_{j,1}(x_1,\dots,x_{j-1})u_1+\mu_{j,2}(x_1,\dots,x_{j-2})u_2+\dots+\\
\vspace{-0.1cm}\\
\hspace{1cm}\mu_{j,j-1}(x_1)u_{j-1}+u_j,\quad j=1,\dots,d-1,\\
\vspace{-0.1cm}\\
w_d=(c_1x_{d-1}+\rho_1(x_1,\dots,x_{d-2}))u_1+(c_2x_{d-2}+\rho_2(x_1,\dots,x_{d-3}))u_2+\\
\vspace{-0.1cm}\\
\hspace{1cm}+\dots+(c_{d-2}x_2+\rho_{d-2}(x_1))u_{d-2}+c_{d-1}x_1u_{d-1},
\end{array}\label{vectorw}
\end{equation}    
where $\mu_{j,k}$, $\rho_i$ are polynomials and $c_{\ell}\in\FF^*$ for all $\ell$. Further, let $v\in\Pi$ and $w':=({\bf 1}+y)v$. Then we have
\begin{equation}
\begin{array}{l}
w'_j=y_{j-1}v_1+y_{j-2}v_2+\dots+y_1v_{j-1}+v_j,\quad j=1,\dots,d-1,\\
\vspace{-0.1cm}\\
w'_d=y_{d-1}v_1+y_{d-2}v_2+y_{d-3}v_3+\dots+y_1v_{d-1}.
\end{array}\label{wprime}
\end{equation}
Equating $w_j$ and $w'_j$ for $j=1,\dots d-1$ in (\ref{vectorw}), (\ref{wprime}), we find
\begin{equation}
\begin{array}{l}
v_j=\eta_{j,1}(y,x_1,\dots,x_{j-1})u_1+\eta_{j,2}(y,x_1,\dots,x_{j-2})u_2+\dots+\\
\vspace{-0.1cm}\\
\hspace{1cm}\eta_{j,j-1}(y,x_1)u_{j-1}+u_j,\quad j=1,\dots,d-1,\\
\end{array}\label{vectorv}
\end{equation}
where $\eta_{j,k}$ are polynomials. Thus, for every $y$ and $u$ we have now chosen $v$ in terms of $x_i$, $i=1,\dots,d-1$. It remains to choose $x_i$ so that $w_d=w'_d$ for all $u\in\Pi$.

Plugging expressions (\ref{vectorv}) into the formula for $w_d'$ in (\ref{wprime}), we obtain
\begin{equation}
\hspace{-0.1cm}\begin{array}{l}
w'_d=y_{d-1}u_1+y_{d-2}\Bigl(\eta_{2,1}(y,x_1)u_1+u_2\Bigr)+\\
\vspace{-0.1cm}\\
\hspace{1cm}y_{d-3}\Bigl(\eta_{3,1}(y,x_1,x_2)u_1+\eta_{3,2}(y,x_1)u_2+u_3\Bigr)+\dots+\\
\vspace{-0.1cm}\\
\hspace{1cm}y_1\Bigl(\eta_{d-1,1}(y,x_1,\dots,x_{d-2})u_1+\eta_{d-1,2}(y,x_1,\dots,x_{d-3})u_2+\dots+\\
\vspace{-0.1cm}\\
\hspace{1cm}\eta_{d-1,d-2}(y,x_1)u_{d-2}+u_{d-1}\Bigr).
\end{array}\label{wprimenew}
\end{equation}
We now compare the coefficients at $u_j$ in formulas for $w_d$ and $w_d'$ in (\ref{vectorw}) and (\ref{wprimenew}). It is straightforward to see that, starting with $u_{d-1}$ and progressing to $u_1$, one can choose $x_i$ in terms of $y$ recursively, starting with $x_1$ and progressing to $x_{d-1}$, in such a way that the coefficients at $u_j$ are pairwise equal for all $j$. Indeed, choose $x_1:=y_1/c_{d-1}\in A_1$, $x_2:=[y_2+y_1\eta_{d-1,d-2}(y,x_1)-\rho_{d-2}(x_1)]/c_{d-2}\in A_2$, etc. This choice of $x_i$, $i=1,\dots,d-1$, guarantees that $w_d=w'_d$ for all $u\in\Pi$. 

We now set $\varphi:=\varphi_{x_{d-1}}\circ\dots\circ\varphi_{x_1}$. Then, for every element $u\in\Pi$ and for the corresponding element $v\in\Pi$ defined by formulas (\ref{vectorv}), we have $\varphi(u)=({\bf 1}+y)v$. Hence $\varphi(\Pi)=({\bf 1}+y)\Pi$ as required.\qed

{\obeylines
\noindent Department of Mathematics
\noindent The Australian National University
\noindent Canberra, ACT 0200
\noindent Australia
\noindent e-mail: alexander.isaev@anu.edu.au
}

\end{document}